\documentstyle[12pt]{article}
\newcommand{\be}{\begin{equation}}
\newcommand{\ee}{\end{equation}}
\newcommand{\ba}{\begin{array}}
\newcommand{\ea}{\end{array}}
\newcommand{\ben}{\begin{enumerate}}
\newcommand{\een}{\end{enumerate}}

\newcommand{\ov}{\overline}
\newcommand{\tr}{{\rm tr}\, }

\newcommand{\Ph}{\mbox{$\bf Ph$}\, }
\newcommand{\im}{\mbox{Im}\, }

\font \msb=msbm10 scaled \magstep1

\newcommand{\bR}{\mbox{\msb R} }
\newcommand{\bC}{\mbox{\msb C} }

\newcommand{\ar}{\alpha }

\newcommand{\gr}{\gamma }

\newcommand{\er}{\varepsilon }

\newcommand{\car}{\ov{\ar}}
\newcommand{\cx}{\ov{x}}

\font \eul=eufm10 scaled \magstep2

\newcommand{\gotG}{\mbox{\eul g}}
\newcommand{\gotH}{\mbox{\eul h}}

\begin{document}

\title{\bf Free motion on the Poisson plane and sphere}
\author{{\bf S. Zakrzewski}  \\
\small{Department of Mathematical Methods in Physics,
University of Warsaw} \\ \small{Ho\.{z}a 74, 00-682 Warsaw, Poland} }

\date{}
\maketitle
\begin{abstract}
Poisson plane and sphere --- homogeneous spaces of Poisson
groups $E(2)$ and $SU(2)$ (resp.) --- have phase spaces
(corresponding symplectic groupoids), in which a free
Hamiltonian is naturally defined.  We solve the equations of
motion and point out some unexpected features: free motion on
the plane is bounded (periodic) and free trajectories on the
sphere are all circles except the big ones.
\end{abstract}

\section{Introduction}

This paper is a continuation of an earlier work
\cite{poi,poican,k-part,zakop,stand} on examples of classical
mechanical systems based on Poisson symmetry. Our examples have
the following common structure: configurations are described by
a Poisson manifold $Q$ which is a quotient of a Poisson group
$G$ by a Poisson subgroup $H$ (we shall assume that both groups
are connected).  In such a case of a Poisson homogeneous
manifold, it is particularly simple to describe its phase space
$\Ph Q$ and the canonical moment map $J\colon \Ph Q\to
\ov{G^*}$, where $G^*$ is the Poisson dual of $G$ (for any
Poisson manifold $P$, we denote by $\ov{P}$ the same manifold with
the opposite Poisson structure). $\Ph Q$ can be identified with
the symplectic reduction of $\Ph G$ with respect to the
constraint
\be\label{constr}
J_H^{\rm right}= \{e\},
\ee
where $J_H^{\rm right}:\Ph G\to H^*$ is the canonical moment for
the (phase lift of the) right translations on $G$ by elements of
$H$, and $e$ is the unit of $H^*$.  Then $J$ is just the
restriction of the canonical moment for the (phase lift of the)
left translations on $G$ to the above constraint ($J$ is
constant on the characteristics because left and right
translations commute).

In favorable cases, $\Ph G = G\Join G^* \cong G\cdot G^*$ is just
the Manin group \cite{qcp} (the group corresponding to the Lie algebra
$\gotG\Join\gotG ^*$ of the Manin triple). The (phase lifts of the)
left and the right translations on $G$ have canonical moment
\be\label{lr}
 J^{\rm left} (g\cdot g^*)= {^g{g^*}}\qquad\mbox{\rm and}\qquad
 J^{\rm right}(g\cdot g^*)=g^*,
\ee
respectively. Here $g\in G$, $g^*\in G^*$,  and $^gg^*\in G^*$ is
the {\em dressing action} of $g$ on $g^*$ defined by
$$ g\cdot g^* = { ^g{g^*}}\cdot g',$$
where $g'\in G$ ($g'$ is denoted by $g^{g^*}$ in \cite{qcp}).
To the (Poisson) inclusion of a Poisson subgroup $H$ in $G$
there corresponds a Poisson projection from $G^*$ to $H^*$. The
kernel of this projection is the subgroup $H^{\circ}$ in $G^*$
corresponding to the annihilator $\gotH ^{\circ}$ of $\gotH$ in
$G^*$.  Composing $J^{\rm right}$ with this projection we obtain
$J_H^{\rm right}$, hence the constraint submanifold in $G\cdot
G^*$ defined by (\ref{constr}) is simply
\be\label{cons}
 K= G\cdot H^{\circ} .
\ee
The phase lift of a right translation by $h\in H$ on $G$ is the
right translation by $h$ on $G\cdot G^*$. Therefore, the
symplectic reduction of $K$ which coincides with the quotient of
$K$ by the action of $H$ is just the set of the right $H$-cosets
in $K$,
\be\label{PhQ}
 \Ph Q = (G\cdot H^{\circ})/H,
\ee
which is a bundle over $Q=G/H$, associated with the principal
$H$-bundle $G$ over $G/H$ and the action of $H$ on $H^{\circ}$ (it is
easy to see that $[\gotH ,\gotH ^{\circ}]\subset \gotH ^{\circ}$
in the Lie algebra of the Manin triple, hence
$h^{-1}H^{\circ}h\subset H^{\circ}$ for $h\in H$ and therefore
$G\cdot H^{\circ}$ is right $H$-invariant). Moreover, if $Q=G/H$
is identified with a submanifold of $G$, one can just identify
\be\label{QH}
 \Ph Q = Q\cdot H^{\circ}.
\ee
Consider now less favorable cases, when we can still assume that
$G$ and $G^*$ can be viewed as subgroups of the Manin group and
$G\cap G^* =\{ e\}$, but $G\cdot G^*$ is not a subgroup. We have
in this case
$$\Ph G = \mbox{connected component of}\;\; \{e\}\;\;
\mbox{in} \;\;\; (G\cdot G^*)\cap (G^*\cdot G)$$
and instead of (\ref{PhQ}) we have
\be\label{XPhQ}
\Ph Q = (\Ph G\cap G\cdot H^{\circ})/H =(G^*\cdot G \cap G\cdot
H^{\circ})/H .
\ee
If $Q=G/H$ is identified with a submanifold of $G$, then
formula (\ref{QH}) is generalized to
\be\label{XQH}
\Ph Q = \Ph G\cap (Q\cdot H^{\circ}).
\ee

The paper is organized as follows. In Section~2 we recall the
results of \cite{poi,poican} concerning the Poincar\'{e}
symmetry in two dimensions.  In Section~3 we consider the
Euclidean symmetry in two dimensions, and in Section~4 we
discuss the case of the ($SU(2)$-symmetric) Poisson 2-sphere.
The above introduced notation is used without further
explanation.

\section{Two-dimensional Minkowski space-time}

In \cite{poi} we have discussed the Poisson symmetry described by
$$G = \{(v^+,v^-,a)\colon \; v^+,v^-,a\in \bR\},\;\;\;
(v^+,v^-,a)(u^+,u^-,b) = (v^+ +e^au^+,v^- + e^{-a}u^-,a+b)$$
$$ H = \{(0,0,a)\colon a\in \bR \},\qquad Q=G/H\cong
\{(x^+,x^-,0)\colon \;  x^+,x^-\in \bR\}$$
($x^+,x^-$ and $v^+,v^-$ are light-cone coordinates), with Poisson brackets
$$ \{v^+,v^-\} =  \er v^+v^-,  \qquad
\{a,v^{\pm} \} =  \er v^{\pm },\qquad
\{x^+,x^-\}=\er x^+x^-.
$$
Then the dual Poisson group is $G^* = \{ (v_+,v_-,b)\colon \;
v_+, v_-,b\in \bR\}$ with multiplication
$$(v_+,v_-,b)(u_+,u_-,b') = (v_+e^{-\er b'/2} +e^{\er
b/2}u_+,v_-e^{-\er b'/2} + e^{\er b/2}u_-,b+b'),$$
and Poisson brackets $\{v_+,v_-\}= 0$, $ \{b,v_\pm \} = \pm
v_\pm $.  We shall denote the elements of $H^{\circ }=\{ (\eta
_+,\eta _-,0)\in G^*\}$ simply by $\eta \equiv (\eta _+,\eta
_-)$ and the elements of $Q$ by $x\equiv (x^+,x^-)$.  For
$(x,\eta )\in\Ph Q$ (formula ({\ref{XQH})), we obtain (using (9)
of \cite{poi}) the groupoid projections
\be\label{gpr}
 (x,\eta )_L= (x^+,x^-),\qquad (x,\eta )_R= \left(
\frac{x^+}{1+\er \eta _-x^-},\frac{x^-}{1-\er \eta _+x^+}\right),
\ee
Poisson brackets (using the fact that the map $(x,\eta)=\xi \mapsto
(\xi _L ,\xi _R)\in Q\times \ov{Q}$ is Poisson)
$$\{x^+,x^-\}=\er x^+x^-,\;\;\; \{\eta _+,\eta _-\}=\er \eta _+\eta _-,
\;\;\; \{ \eta _+, x^+\} = 1 - \er \eta _+x^+,$$
$$\{\eta _+,x^-\}=-\er \eta _+x^- ,\;\;\;
\{\eta _-,x^+\}=\er \eta _-x^+,\;\;\;
 \{ \eta _-, x^-\} = 1 +\er \eta _-x^-
 $$
and the moment map $J\colon \Ph Q\to \ov{G^*}$ (using again (9)
of \cite{poi} and the change of coordinates in \cite{poican})
which associates with $(x,\eta )\in \Ph Q$ the following element 
of $G^*$: 
$$ \left( \eta _+\sqrt{\frac{1+\er \eta _-x^-}{1-\er \eta _+x^+}},
\eta _-\sqrt{\frac{1+\er \eta _-x^-}{1-\er \eta _+x^+}},
-\frac1{\er }\log (1+\er \eta _-x^-)(1-\er \eta _+x^+)\right) .
$$
Transporting the Casimir function $v_+v_-$ on $G^*$ to $\Ph Q$
by $J$ we obtain the mass shell in the usual form $\eta _+\eta _-=m^2$.
Calculating the Poisson brackets of coordinates with the `Hamiltonian'
equal $\eta _+\eta _-$ yields the following equations of motion of
a free particle
$$
\dot{x}^+  =  \eta _-,\;\;\;
\dot{x}^-  =  \eta _+,  \;\;\;
\dot{\eta}_+  = - \er  \eta _-\eta _+^2,  \;\;\;
\dot{\eta}_-  =  \er  \eta _+\eta _-^2 .
$$
It follows that the world line of the particle is a hyperbole
\be\label{hip}
 (x^+ - c^+)(x^- - c^-)= - \frac1{\er ^2 m^2}.
\ee

Note that in $\Ph Q$ one can introduce commuting space-time
coordinates as follows. Since the Poisson bivector on $Q$
is built of commuting vector fields:
\be\label{commut}
\pi = \er x^+\frac{\partial}{\partial x^+}\wedge x^-
\frac{\partial}{\partial x^-},
\ee
one can realize (see \cite{poican,abel})  $\Ph Q$
in the cotangent bundle $T^*Q$. For the usual coordinates
 $(q,p)=(q^+,q^-,p_+,p_-)$ of $T^*Q$ we have the usual canonical
commutation relations, the groupoid projections from $\Ph Q$ to
$Q$ being given by 
\be\label{cgpr}
(q,p)_L=(q^+e^{\frac{\er}{2}p_-q^-},q^-e^{-\frac{\er}{2}p_+q^+}),\qquad
(q,p)_R=(q^+e^{-\frac{\er}{2}p_-q^-},q^-e^{\frac{\er}{2}p_+q^+}).
\ee
Identifying projections (\ref{gpr}) and (\ref{cgpr}) we obtain
formulae relating $(x,\eta )$ and $(q,p)$. In particular, we obtain
expressions for the moment map and the mass shell in the new variables:
$$
J(q^+,q^-,p_+,p_-) = ( P_+,P_-,\Pi _+ - \Pi _-),\qquad  P_+P_-=m^2,
$$
where
$$
P_+ = \frac{\sinh \frac{\er}{2} p_+q^+}{\frac{\er}{2} q^+ },\qquad
P_- =\frac{\sinh\frac{\er}{2} p_-q^-}{\frac{\er}{2} q^- },\qquad
\Pi _+=p_+q^+,\qquad \Pi _-=p_-q^-.
$$
This gives the  equations of motion
$\dot{P}_{\pm}=0$, $ \dot{\Pi} _{\pm}=m^2$,
which yield now world lines different from hyperboles (\ref{hip}):
$$
q^+ = e^{a}\, \frac{\sinh \frac{\er}{2} \tau }{\frac{\er}{2} m },\qquad
q^- =e^{-a}\, \frac{\sinh \frac{\er}{2} (\tau - b )}{\frac{\er}{2} m },
$$
where $\tau$ is the parameter and $a$, $b$ are some constants.

Note that the cotangent bundle projection turns out to be a kind of
geometric mean of the left and right groupoid projections:
$$ x_L^+x_R^+=(q^+)^2,\qquad x_L^-x_R^-=(q^-)^2 .$$

\section{The Poisson plane}

We identify the (double cover of the) Euclidean group $E(2)$
in two dimensions with the set of matrices
$$G=\left\{\left(\ba{cc} \ar & 0 \\ \gr & \car \ea\right)\;
: \;\;\;\ar ,\gr \in \bC,\;\; \ar \car =1 \right\} .$$
The standard Poisson structure ($\er $ is the deformation parameter) 
\be\label{E2}
 \{ \ar , \gr \} = i\er \ar\gr,\qquad \{ \gr,\ov{\gr}\} =0 ,
\ee
corresponds to the Manin group $(SL(2,\bC ); G, G^* )$ with 
the scalar product defined by
\be\label{scap}
 \frac{1}{\er} \im \tr XY  \qquad \mbox{for}\;\; X,Y\in sl(2,\bC),
\ee
and the dual Poisson group being
\be\label{dual}
G^* = \left\{\left( \ba{cc} \rho & n \\ 0 & \rho ^{-1}\ea\right)\;
 : \;\;\;\rho >0,\; n\in \bC \right\} .
\ee
The Poisson brackets on $G^*$ are then \ $\{ \rho ,n\} = -i\er
\rho n $, $ \{ n,\ov{n} \} = 0$. \ 
It is convenient to use new parameters $(P,s)$ on  $G^*$, given by
$$  \rho = \exp (\er s),\qquad  n = i\er \ov{P}.$$
In the limit $\er \to 0$, these parameters become the usual translational
and rotational momenta, respectively. Their Poisson brackets are
simply
$$ \{ s ,P\} = i\er P,\qquad \{ P,\ov{P} \} = 0.$$

Dividing $G$ by the Poisson subgroup $H$ composed of diagonal
matrices, we get the quotient $Q=G/H$ which may be identified with
the complex plane:
$$ Q\cong \left\{ \left( \ba{cc} 1 & 0\\ x & 1\ea\right) \; : \;\;\;
x\in \bC \right\} .$$
Using (\ref{E2}) it is easy to find the Poisson structure induced
on the plane by the projection from $G$ on $Q$ (given by $x=\car \gr$):
\be\label{cxx}
\{ \cx ,x\} = 2i\er |x|^2 .
\ee
If $x=x^1+i x^2$, then $\{ x^1,x^2\} = \er ((x^1)^2 + (x^2)^2)$.

The subgroup $H^{\circ}$ in $G^*$ is composed of elements of the
form $(P,s)=(\eta,0)$, where $\eta\in \bC $. We shall denote
$(\eta ,0)$ by $\eta$. Writing
$$ \left( \ba{cc} 1 & 0\\ x & 1\ea\right)
\left( \ba{cc} 1 & i\er \ov{\eta} \\ 0 & 1\ea\right)$$
as a product $g^* g$, where $g\in G$, $g^*\in G^*$, we  get the
groupoid projections
\be\label{gpr1}
 (x,\eta )_L=x,\qquad (x,\eta )_R=\frac{x}{1-i\er \cx\eta},
\ee
defined on $\Ph Q$ as given in (\ref{XQH}), and the moment map
\be\label{JE2}
J(x,\eta ) = \left( \eta\cdot \frac{|1-i\er\cx\eta |}{1-i\er\cx\eta },
-\frac{1}{\er }\log |1-i\er\cx\eta |\right) .
\ee
Using the fact that $\Ph Q\ni \xi\mapsto (\xi _L,\xi _R)\in
Q\times  \ov{Q}$ is a Poisson map, one gets the Poisson
structure of $\Ph Q$: in addition to (\ref{cxx}) we have
\be\label{phQ}
\{ \ov{\eta},\eta \}=-2i\er |\eta |^2,\qquad \{ \eta ,x\} =
-2i\er\eta x,\qquad \{\eta ,\cx \} = 2(1-i\er \cx \eta ).
\ee
Now we note that the Casimir function $|P|^2$ on $G^*$ equals
$|\eta |^2$, when transported to $\Ph Q$ by the moment map.
Therefore it is natural to consider the Hamiltonian
$$ {\cal H} = \frac12 |\eta |^2 $$
as describing the (`$\er$-analogue' of the) free dynamics.
Integrating the equations of motion
\be\label{eom}
\dot{x} = \{ {\cal H},x\} = \eta,\qquad \dot{\eta}=-i\er |\eta |^2\eta,
\ee
we obtain easily
$$\eta = \eta _0 \exp (-2i\er Et),\qquad
x= x_0 + \eta _0 \frac{1-\exp (-2i\er Et)}{2i\er E}, \qquad
E\equiv {\cal H}=\frac12 |\eta _0|^2.
$$
We now see the `effect of the deformation': free trajectories are
bounded. They are circles or points (never straight lines). The radius
of the circle is inversely proportional to the velocity!

As in Section~1, one can introduce (somehow externally)
commuting positions in $\Ph Q$, using the following
representation of the Poisson bivector on $Q$:
$$ \pi =\er |x|^2 \partial _1\wedge \partial _2 =
\er (x^1\partial _1 +x^2\partial _2)\wedge (x^1\partial _2
-x^2\partial _1), \;\;\;\;\mbox{where}\;\;
\partial _j:=\frac{\partial}{\partial x^j},\; j=1,2.
$$
For $(q,p)=(q^1+iq^2,p_1+ip_2)\in T^*Q$ we have the usual
canonical commutation relations (the only non-zero are
$\{p,\ov{q}\}=2$) and we get the following groupoid projections
\be\label{cgpr1}
 (q,p)_L=\exp (-i\frac{\er}{2} \ov{q}p)\cdot q,\qquad
(q,p)_R=\exp (i\frac{\er}{2} \ov{q}p)\cdot q .
\ee
Identifying projections (\ref{gpr1}) and (\ref{cgpr1}) we obtain
formulae relating $(x,\eta )$ and $(q,p)$. In particular, we obtain
expressions for the moment map and the Hamiltonian in the new variables:
$$
J(q,p) = ( P, -\im \ov{q}p),\qquad  {\cal H} = \frac12 P\ov{P},
\qquad \mbox{where}\;\;
P = \frac{\sin \frac{\er}{2} \ov{q}p}{\frac{\er}{2} \ov{q}} .
$$
Since $\{P,\ov{P}\}=0$, the `effective' momentum $P$ is
conserved. Using $\{P,\ov{q}p\}=2P$ and $\{q\ov{p},\ov{q}p\}
=0$, we get $ \frac{d}{dt}(\ov{q}p) = |P|^2 = 2E$ and
$$
q=\frac{\sin \frac{\er}{2} q\ov{p}}{\frac{\er}{2} \ov{P}}=
q_0 \frac{\sin \frac{\er}{2} ((q\ov{p})_0 + 2Et)}{\sin
\frac{\er}{2} (q\ov{p})_0} .$$
Trajectories in terms of commuting positions $q$ are not circles
but more complicated closed curves. They arise as geometric mean
of two circular motions, since
$$ \xi _L\cdot \xi _R =(q,p)_L\cdot (q,p)_R = q^2 \qquad\mbox{for}
\;\; \xi\in \Ph Q.$$
Using this property, we obtain also a remarkable formula for the
Hamiltonian:
$$ 2\er ^2{\cal H} = |\er P|^2 
= \left| \frac{2\sin \frac{\er}{2} \ov{q}p}{q}\right|^2
= \left| \frac{ \exp (i \frac{\er}{2} \ov{q}p) q  -
 \exp (-i \frac{\er}{2} \ov{q}p) q }{q^2}\right|^2 $$
$$
= \left| \frac{\xi _R - \xi _L}{\xi _L\xi _R}\right|^2 =
\left| \frac{1}{\xi _L} - \frac{1}{\xi _R}\right|^2 .
$$

\section{Free motion on the Poisson sphere}

We start with simple observations relating the usual free motion
on the sphere $Q=G/H=SU(2)/{S^1}$ to the free motion on $G=SU(2)$.
The Hamiltonian $H\colon T^*G\to \bR $ of the latter may be written
as
\be\label{free}
{\cal H}(\xi ) = \frac12  J^2,
\ee
where $J: T^*G\to \gotG ^*$ is the moment map for the right
translations ($J(\xi ) = g^{-1}\xi $ for $\xi\in T_g^*G$) and
$J^2(\xi )$ denotes the scalar square of $J(\xi )$ in terms of
some invariant scalar product on $\gotG ^*$. Specifically, using
the natural orthogonal basis
$$
J_1=\left(\ba{cc} 0 & i \\ i & 0\ea\right),\qquad
  J_2 =\left(\ba{cc} 0 & -1 \\ 1 & 0\ea\right),\qquad
  J_3 = \left(\ba{cc} i & 0 \\ 0 & -i\ea\right),\qquad
$$
of $su(2)$, denoting by $\tilde{J}_k$ the corresponding
left-invariant vector fields on $G$ and (again) by $J_k$ the
associated functions on $T^*G$,
$$ J_k (\xi ) = \left\langle \xi , \tilde{J}_k (g)\right\rangle
\qquad \mbox{for}\;\; \xi\in T_g^*G ,$$
we may set
$$ {\cal H} = \frac12 (J_1^2 + J_2^2 + J_3^2).$$
It is well known that the resulting free trajectories on
$SU(2)\cong S^3$ are `big circles' (translated one-parameter
subgroups).

Let $H\cong S^1$ be the subgroup generated by $J_3$. To the
reduction map $G\to Q=G/H $ there corresponds the symplectic reduction
\be\label{redu}
 T^*G\supset K \to T^*Q ,
\ee
where the constraint set $K$ is given by $\{ J_3 =0\}$. Since
the $S^1$ action from the right (generated by $J_3$) preserves
the Hamiltonian (i.e. $\{J_3,{\cal H}\}=0$), projections of
trajectories of $H$ are trajectories of the projected
Hamiltonian ${\cal H}_{\rm red}$ on $T^*Q$. Moreover (since
reduction (\ref{redu}) is implied by the configurational
reduction), also configurational trajectories on $SU(2)$ are
projected on configurational trajectories on $Q=S^2$.

The following lemma is easily proved.

\vspace{2mm}

\noindent
{\bf Lemma.} \ {\em Projecting all `big circles' from $SU(2)$ to
$S^2$, one obtains all circles on $S^2$.}

\vspace{2mm}

It is perfectly known that free trajectories on $S^2$ are not
`all' circles but only the `big' ones (or points --- in the case
of rest). This may seem to
contradict the statement before the lemma. There is no paradox,
of course, because we project from $T^*G$ only those
trajectories which are in the constraint $K$.
These are exactly those trajectories whose configurational
velocity is perpendicular to $\tilde{J}_3$ (the direction of the
$S^1$ action):
$$ J_3 (\xi )=0\Longleftrightarrow
\left\langle \xi , \tilde{J}_3 (g)\right\rangle
=0 \Longleftrightarrow v \perp \tilde{J}_3$$
($v\in T_gG$ is obtained from $\xi\in T_g^*G$ by using the
invariant metric on $G$). One can easily check that this
perpendicularity condition is exactly equivalent to the fact that the
projection on $S^2$ is a `big' circle (or a point).

Now we consider the standard Poisson structure on $G=SU(2)$,
corresponding to the Manin group $(SL(2,\bC );G,G^*)$, with the
scalar product and the dual group as before, given by (\ref{scap})
and (\ref{dual}).
In \cite{zakop} we have calculated the symplectic structure of
$\Ph G =SL(2,\bC )$ and solved the equations of motion for the
following analogue of the free Hamiltonian (\ref{free}):
$$ {\cal H}(\xi ) = \frac12 \tr \xi ^{\dagger}\xi ,\qquad \xi\in
SL(2,\bC).$$
We shall write the result for the rescaled Hamiltonian
\be\label{free1}
 \tilde{\cal H} := \frac{{\cal H} -1}{4\er ^2},
\ee
because this one tends to (\ref{free}) when $\er \to 0$,
while using the following parame\-tri\-za\-tion of $G^*$:
$$ G^*\ni b=\left( \ba{cc} \rho & n \\ 0 & \rho ^{-1}\ea\right) =
\left( \ba{cc} \exp (\er s) & 2\er w \\ 0 & \exp (-\er s)
\ea\right) .
$$
The result is that the phase trajectory $t\mapsto \xi = gb\in
G\cdot G^*=SL(2,\bC )$ satisfies
$$g^{-1}\dot{g} = {\cal F} (b):=\frac{i}{2}
\left( \ba{cc} (2\er )^{-1}\sinh 2\er s \, +\er |w|^2 &
 w\exp (-\er s) \\ \ov{w}\exp (-\er s) &
 - (2\er )^{-1}\sinh 2\er s \, -\er |w|^2 \ea\right)
$$
and $\dot{b}=0$.
This means that the left groupoid projection $t\mapsto g(t)$ of
the phase trajectory is the usual `big' circle, and the motion
is uniform with constant velocity given by the deformed
`Legendre transformation' $G^*\ni b\mapsto v:={\cal F}(b)\in \gotG $.

The free Hamiltonian (\ref{free1}) is a Casimir function on
$G^*$ and commutes with left and right momenta (\ref{lr}). In
particular, it commutes with the constraint (\ref{constr}), and,
consequently, may be projected down to the reduced space $ \Ph Q
= \Ph (G/H) $. It is this function which we consider as the
$\er$-analogue of the Hamiltonian of the free motion on $S^2$.
Note that the resulting trajectories on $S^2$ being left
projection of trajectories in $\Ph Q$, are at the same time the
projections of those `big' circles on $SU(2)$ which come from
phase trajectories living in $K =G\cdot H^{\circ}$ (formula
(\ref{cons})). This follows from the fundamental property of
morphisms of groupoids \cite{qcp} (morphisms commute with
groupoid projections).  Phase trajectory lives in $K$ if and
only if $b\in H^{\circ}$, i.e.
$$ b = \left( \ba{cc} 1 & 2\er w \\ 0 & 1 \ea\right) .$$
For such $b$,
$$ {\cal F} (b)=\frac{i}{2}
\left( \ba{cc} \er |w|^2 &  w \\
 \ov{w} & -\er |w|^2 \ea\right)
$$
has vanishing scalar product with $J_3$ only when $w=0$. It
follows that the velocity is perpendicular to the $S^1$ action
only if it is zero. This means that we never get a `big' circle
on $S^2$. Moreover, since velocities may have here any angle
(different from $0$ and $\pi /2 $) with $J_3$, we get all circles
on $S^2$, except the big ones. This fact corresponds to the
previous result concerning the Poisson plane: trajectories could
be any circles `except' the straight lines.


\begin{thebibliography}{99}

\bibitem{poi} S. Zakrzewski,
``Poisson space-time symmetry and corresponding elementary
systems'', in: {\em Quantum Symmetries}, Proceedings of the II
International Wigner Symposium, Goslar 1991, H.D.~Doebner and
V.K.~Dobrev (Eds.), pp. 111--123.
\bibitem{poican} S. Zakrzewski, ``Poisson Poincar\'{e} particle
and canonical variables'', in: {\em Generalized Symmetries},
Proceedings of the International Symposium on Mathematical
Physics, Clausthal, July 27--29, 1993, H.-D. Doebner, V.K.
Dobrev and A.G. Ushveridze (Eds.), 1994, pp. 165--171.
\bibitem{k-part} S.~Zakrzewski, ``On the classical $\kappa
$-particle'',
in: {\em Quantum Groups, Formalism and Applications},
Proceedings of the XXX Winter School on Theoretical Physics
14--26 February 1994, Karpacz, J.~Lukierski, Z.~Popowicz,
J.~Sobczyk (eds.), Polish Scientific Publishers PWN, Warsaw
1995, pp.~573--577. Also: hep-th{/}9412098.
\bibitem{zakop} S.~Zakrzewski, ``Classical mechanical systems
based on Poisson symmetry'', in: Proceedings of the
II German-Polish Symposium, Zakopane, September 1995,
{\sl Acta Physica Polonica B}, vol. {\bf 27} No.~10 (1996),
2801--2810.
\bibitem{stand} S. Zakrzewski, ``Phase spaces related to standard
classical $r$-matrices'', {\em J.~Phys.~A}: Math. Gen. {\bf 29} 
(1996) 1841--1857.
\bibitem{qcp}  S. Zakrzewski, ``Quantum and classical
pseudogroups. Part I and II'',
{\em Comm. Math. Phys.} {\bf 134} (1990), 347--395.
\bibitem{abel} S. Zakrzewski, ``Geometric quantization of
Poisson groups --- diagonal and soft deformations'',
Proceedings of the Taniguchi Symposium {\em
Symplectic geometry and quantization problems}, Sanda (1993),
Y.~Maeda, H.~Omori and A.~Weinstein (Eds.), Contemporary
Mathematics {\bf 179}, 1994, 271--285.



\end{thebibliography}
\end{document}